\documentclass[12pt]{amsart}
\usepackage{amscd}
\usepackage{amssymb}
\theoremstyle{plain}
\newtheorem{theorem}{Theorem}
\newtheorem{proposition}{Proposition}[section]
\newtheorem{algorithm}{Algorithm}

\newtheorem{corollary}[proposition]{Corollary}

\newtheorem{lemma}[proposition]{Lemma}

\newcommand{\TeXButton}[1]{{}}
\numberwithin{equation}{section}

\newcounter{stepletter}    
\newenvironment{steplist}{\begin{list}{{\sc Step} \arabic{stepletter}:}
{\usecounter{stepletter}\setlength{\labelsep}{1em}\setlength{\labelwidth}{7em}%
\setlength{\leftmargin}{5em}}}{\end{list}} 

\newcounter{typeletter}   \renewcommand{\thetypeletter}%
{{\sc Type}~\arabic{typeletter}}

\newcounter{lcletter}

\newcounter{romnum}

\begin{document}
\author{Barbara L. Osofsky }
\address{Department of Mathematics \\
Rutgers, The State University of New Jersey \\
110 Frelinghuysen Road \\
Piscataway, NJ 08854-8019}
\email{osofsky@math.rutgers.edu}
\dedicatory{This paper is dedicated to the memory of Richard Pierce.}
\date{Uploaded July 14, 2000}
\title{Projective Dimension is a Lattice Invariant}
\subjclass{Primary 13D05, 20K99; Secondary 06E20}
\keywords{Algorithms for free abelian groups, Lifting direct sum decompositions, Projective dimension, commuting idempotents, lattice invariants}

\begin{abstract}
We show that, for a free abelian group $G$ and prime power $p^\nu$, 
every direct sum decomposition of the group $G\left/p^\nu G\right.$
lifts to a direct sum decomposition of $G$.  This is
the key result we use to show that, for $R$ a commutative
von Neumann regular ring, and $\mathcal{E}$ a set of idempotents in
$R$, then the projective dimension of the ideal $\mathcal{E}R$ as an 
$R$-module the same as the projective dimension
of the ideal $\mathcal{E}\mathcal{B}$ as a $\mathcal{B}$-module, where
$\mathcal{B}$ is the boolean algebra generated by $\mathcal{E}\cup\left\{1\right\}$.
This answers a thirty year old open question of R. Wiegand. 
\end{abstract}

\maketitle

\section{Introduction.}

Back in the late 60's, Roger Wiegand asked the following question in \cite
{wiegand}:

\begin{quotation}
{Let $R$ be a commutative [von Neumann] regular ring and $J$ an ideal of $R$ generated by
a set $\mathcal{E}$ of idempotents. Let $\mathcal{B}$ be the Boolean algebra
of all idempotents of $R$. Then is the projective dimension of $J=
\mathcal{E}R$ as an $R$-module the same as the projective dimension of $
\mathcal{EB}$ as a $\mathcal{B}$-module? }
\end{quotation}

In this paper we show that the answer to this question is `yes'.

Richard Pierce popularized this problem, and did some of the early
work on it. It is not difficult to see that the answer is `yes' if $J$
is projective. In \cite{pierce:67} Pierce showed that projective
dimension of an ideal generated by an independent set of idempotents
in a boolean algebra was $\kappa$ where the independent set had
cardinality $\aleph_\kappa$ (here $\kappa \ge \omega$ is replaced by
$\infty$ for projective dimension).  Osofsky
\cite{pddirprod} proved the same result for arbitrary commuting
idempotents in any ring, so in the case of ideals generated by
independent idempotents the answer to the Wiegand question is
`yes'. Then Richard Pierce \cite{pierce:76} showed that it is `yes' in
case either the projective dimension of $\mathcal{E}R$ or the projective 
dimension of $\mathcal{EB}$ is one. Since then, the problem has
been solved in some special cases with extra hypotheses on the
idempotents forcing projective dimension to be the subscript of the minimal
aleph of a generating set, although the general problem remained open. 

The essence
of the problem is that the additive order of some of the idempotents
in $\mathcal{E}$ might be one prime (for example the prime 2 in case
$R=\mathcal{B}$) and a different prime in another ring $ R^{\prime}$,
or perhaps even infinite in a third ring $R^{\prime\prime}$.  Here we
conquer the problem of different primes by working in a (not regular)
ring $R$ of characteristic 0.
We show that the answer to Wiegand's question is `yes'
in all cases. 

In Section~\ref{Part2}, we prove a subtle
but elementary result about free abelian groups, namely for any free abelian group
$G$ and any direct sum decomposition of $G\left/p^\nu G\right.$, this decomposition lifts
to a direct sum decomposition of $G$. In Section~\ref{Part3} we apply 
 this result to any commutative von Neumann regular ring $R$ containing a lattice
of idempotents isomorphic to $\mathcal{B}$. Unlike Pierce's paper
concerning the case of projective dimension 1 (\cite{pierce:76}), we
do not give an internal characterization of projective dimension of
ideals in a commutative von Neumann regular ring.  However, there is a
candidate for such a characterization in a series of papers by the
author: \cite{ohio:96}, \cite{murcia}, and \cite {ohio:99}.

\section{A theorem on abelian groups}\label{Part2}

The aim in this section is to prove the following: 

\begin{theorem} \label{mainthm} Let $G$ be a free abelian group and $\left\{\overline{b_\alpha}:\alpha\in\mathfrak{I}\right\}$ a (free) basis for $G\left/p^\nu G\right.$ with
$p$ a prime.
Then there exists a family of integers
$\left\{{u_\alpha}:\alpha\in\mathfrak{I}\right\}$, relatively prime to
$p$, and a 
free basis of $G$, $\left\{{y_\alpha}:\alpha\in\mathfrak{I}\right\}$, such that 
$\overline{y_\alpha} =
u_\alpha\,\overline{b_\alpha}$ in $G\left/ p^\nu G\right.$ 
for all $\alpha\in\mathfrak{I}$.
\end{theorem} 

A way of restating this theorem is that the direct sum decomposition $G\left/p^\nu G\right.=\bigoplus_\alpha\, \overline{b_\alpha}\, \mathbb{Z}\left/p^\nu\mathbb{Z}\right.$
lifts to a direct sum decomposition $G=\bigoplus_\alpha y_\alpha\,\mathbb{Z}$.
In fact, any direct sum decomposition of $G\left/p^\nu G\right.$ will lift to a
direct sum decomposition of $G$ by taking bases of each of the summands and lifting
them. We use the fact that the ring $\mathbb{Z}\left/p^\nu\mathbb{Z}\right.$ is local, 
that is, has a unique maximal ideal. If $p^\nu$ is replaced by an arbitrary integer 
which has at least two distinct prime factors, the result is false since $\mathbb{Z}$ is indecomposable 
whereas  $\mathbb{Z}\left/n\mathbb{Z}\right.$ decomposes if $n$ is a product of two relatively 
prime factors $>1$.

\subsection{Basic notation.}

We fix a prime power $p^\nu$. For any abelian group $G$, 
we denote the natural map from $G$ to $G\left/p^\nu G\right.$
by an overline. If $\overline{x}$ is an element of $\overline{G}=G\left/p^\nu G\right.$
we will assume from the notation that $x\in G$ is some preimage of $\overline{x}$.
If $G$ is some free abelian group, we will denote some free basis for $G$ by 
\[ \mathfrak{X}=\left\{x_\sigma:\sigma\in\mathfrak{K}\right\}.\]  and we will denote a basis of 
$\overline{G}$ as a (free) $\overline{\mathbb{Z}}$-module by 
\[ \mathfrak{B}=\left\{\overline{b_\alpha}:\alpha\in\mathfrak{I}\right\}. \]

\subsection{Reduction to the countable case.\label{kap}}

Much of this paper relies heavily on a beautiful 
paper by Kaplansky (\cite{Kaplansky}) for both technique and results. Here we
adapt the basic technique of Kaplansky's paper to get a specialized result on free
abelian groups.  We have the same objective as Kaplansky did, namely to reduce
the question under study to the countable case.

\begin{lemma} \label{OneSummand} Let $G$ be a nonzero free abelian group
with free basis $\mathfrak{X}$, and let $\mathfrak{B}$
be a basis of $\overline{G}$ as a (free) $\overline{\mathbb{Z}}$-module.  
Let $\mathfrak{c}$ be any countable subset of $\mathfrak{B}$. Then there
exists a nonzero countably generated direct summand $H$ of $G$ such that 
\[ \overline{H}=\sum_{i=0}^\infty\,\overline{b_{\alpha_i}}\,\overline{\mathbb{Z}} \]
for $\left\{\overline{b_{\alpha_i}}:i\in\omega\right\}$
some countable subset of $\mathfrak{B}$ containing $\mathfrak{c}$.  
Moreover, $H$ itself is generated by a countable subset of $\mathfrak{X}$.
\end{lemma}

\begin{proof}   
We are given that $\mathfrak{X}=\left\{x_\sigma:\sigma\in\mathfrak{K}\right\}$ is a free basis for $G$. Fix a lifting $\left\{ b_\alpha\right\}$ of $\mathcal{B}$. For any countable subset
$\mathfrak{c}\subseteq \mathfrak{I}$, let $X_\mathfrak{c}\subseteq \mathfrak{K}$ be the smallest (necessarily countable) subset of $\mathfrak{K}$ such that $\sum_{\alpha\in\mathfrak{c}}\,b_\alpha \mathbb{Z}
\subseteq \sum_{\sigma\in X_{\mathfrak{c}}}\,x_\sigma \mathbb{Z}$.  Similarly, for any countable
subset  $\mathfrak{c}'\subseteq \mathfrak{K}$, let $B_{\mathfrak{c}'}\subseteq \mathfrak{I}$ be the smallest
(necessarily countable) subset of $\mathfrak{I}$ such that $\sum_{\sigma\in\mathfrak{c}'}\,{x_\sigma}\, {\mathbb{Z}} 
\subseteq \sum_{\alpha\in B_{\mathfrak{c}'}}\,{b_\alpha} \,{\mathbb{Z}}$.  

Now start with any nonempty countable set $\mathfrak{c}_0$ such that $\mathfrak{c}\subseteq\mathfrak{c_0}\subseteq\mathfrak{I}$. We use finite induction to define two
sequences $\left\{\mathfrak{c}_i,\,\mathfrak{c}_i': i<\omega\right\}$ of countable sets by
\begin{eqnarray*}
\mathfrak{c}_n'&=&X_{\mathfrak{c_n}}     \\
\mathfrak{c}_{n+1}&=&B_{\mathfrak{c}_n'}.  
\end{eqnarray*}
In words, think of $\mathfrak{B}$ as images of  $\left\{{b_\alpha}:\alpha\in\mathfrak{I}\right\}$.   
Starting with a countable subset $\mathfrak{c}_0$ of the basis $\mathfrak{B}$ of $\overline{G}$, use our lifting of $\mathcal{B}$ to get an inverse image $c_0\subseteq G$ and take the smallest countable subset $\mathfrak{c}_0'$
of the basis $\mathfrak{X}$ of $G$ whose span contains $c_0$. Now take images of 
$\mathfrak{c}_0'$ modulo $p^\nu$ and find the smallest countable subset $\mathfrak{c_1}\supseteq\mathfrak{c}_0$ of the basis $\mathfrak{B}$ which span a group containing all 
of the elements of $\overline{\mathfrak{c}_i}$. Iterate a countable number of times. 

We then have  for all $i$, $\mathfrak{c}'_i \subseteq \mathfrak{c}'_{i+1}$ and 
\begin{equation} \tag{$*$}
\overline{F_n}=\sum_{\alpha\in\mathfrak{c}_n}\, \overline{b_\alpha}\,\overline{\mathbb{Z}}\quad
\subseteq \quad \overline{G_n}=\sum_{\sigma\in\mathfrak{c}_n'}\,\overline{x_\sigma}\,
\overline{\mathbb{Z}}\quad \subseteq \quad \overline{F_{n+1}}=\sum_{\alpha\in\mathfrak{c}_{n+1}}\, \overline{b_\alpha}\,\overline{\mathbb{Z} }.
\end{equation}  

Set $H= \sum_{\sigma\in\bigcup_{n=0}^\infty\,\mathfrak{c}_n'} x_\sigma\,\mathbb{Z}$.  Clearly 
$H$ is a direct summand of $G$.  Moreover, $H$ is countably generated since the
indexing set is a countable union of countable sets. Equation~($*$) forces 
\[ \overline{H}=\sum_{\alpha\in\bigcup_{i=0}^\infty\,F_i}\,
\overline{b_\alpha}\,\overline{\mathbb{Z}}.  \]  \end{proof}
% \enlargethispage{\baselineskip}

\begin{lemma} \label{appl} Let $G$ be a nonzero free abelian group, 
and let $\mathfrak{B}$
be a basis of $\overline{G}$ as a (free) $\overline{\mathbb{Z}}$-module.  Then $G$ is the union
of a well-ordered (by inclusion) family $\left\{H_\mu: \mu<\Omega\right\}$ of subgroups
such that:  $H_\mu$ and  $\bigcup_{\kappa<\mu}\, H_\kappa$ are
direct summands of $G$ for every $\mu$ in the ordinal $\Omega$;   for each
$\mu$, 
$H_\mu\left/\bigcup_{\kappa<\mu}\, H_\kappa\right.$ is countable;  
 and each $\overline{H_\mu}$ is generated by some subset of the  
$\left\{\overline{b_\alpha}:\alpha\in\mathfrak{I}\right\}$.\end{lemma}

\begin{proof}  Fix a basis $\mathfrak{X}$ of $G$.
 Well order $\mathfrak{B}$.  Assume
we have $H_\kappa$ for all $\kappa<\mu$ such that:

\begin{enumerate}\renewcommand{\labelenumi}{({\it \roman{enumi}})}
\item Each $H_\kappa$ is
generated by a subset of $\mathfrak{X}$;
\item 
 $\overline{H_\kappa}$ is generated by some subset of $\mathfrak{B}$; and 
\item 
$H_\kappa\supset H_{\kappa'}$ if $\kappa>\kappa'$.  
\item $H_\kappa\left/\bigcup_{\kappa'<\kappa}H_{\kappa'}\right.$ is countably generated.
\end{enumerate}

$H_\kappa$ and $\bigcup_{\mu<\kappa}\,H_\mu$ are direct summands of $G$ since they are generated by  subsets of our fixed basis. 
If $\bigcup_{\kappa<\mu}\, H_\kappa\ne G$, that union cannot map 
 onto $\overline{G}$.  Let $\overline{b_\beta}$ be the smallest element of $\mathfrak{B}$ 
(under the well ordering of $\mathfrak{B}$) not in $\overline{\bigcup_{\kappa<\mu}\,
H_\kappa}$.  Apply Lemma~\ref{OneSummand} to get a countably generated
subgroup $K_\mu$ generated by elements of
$\mathfrak{X}$ with $\overline{b_\beta}\in\overline{K_\mu}$ and 
$\overline{K_\mu}$ generated by a subset of  $\mathfrak{B}$.
Set $H_\mu=K_\mu+\bigcup_{\nu<\mu}\,
H_\nu$.  Since $H_\mu$ clearly has the required properties and this process must eventually
give all of $G$ (at least by the order type of $\mathfrak{B}$), 
by transfinite induction we are done. 
\end{proof}

\begin{corollary} \label{cgenuf} Assume that, for any countably generated free abelian group $G$ with $\mathfrak
{B}$ a basis for $\overline{G}$, there is a direct decomposition lifting of
\[ \overline{G}=\bigoplus_{\alpha\in\mathfrak{I}}\,\overline{b_\alpha}
\,\overline{\mathbb{Z}}
\]
to the direct decomposition
\[ G= \bigoplus_{\alpha\in\mathfrak{I}}\,y_\alpha\mathbb{Z}.
\]
 Then Theorem~\ref{mainthm} is true for any free abelian group $G$.
\end{corollary}  

\begin{proof}  
Using the notation of  Lemma~\ref{appl}, we let
 $G=\bigcup_{\mu<\Omega}\, H_\mu$ where for all $\mu<\Omega$,
$H_\mu=K_\mu +\bigcup_{\kappa<
\mu}\, H_{\kappa}$ with $K_\mu$ countably generated.  
For each $b_{\alpha_i}$ in $K_\mu \smallsetminus \bigcup_{\kappa<
\mu}\, H_{\kappa}$, set $b_{\alpha_i} = c_i + b_i'$, where $c_i$ is
the projection of $b_{\alpha_i}$ to $\bigcup_{\kappa<
\mu}\, H_{\kappa}$.  If $b'_i=0$, ignore it and renumber.
 By assumption, we can lift the 
direct sum decomposition of the quotient
\[ \overline{H_\mu\left/\bigcup_{\kappa<
\mu}\, H_{\kappa}\right.}\approx \overline{K_\mu \left/ K_\mu \cap \bigcup_{\kappa<
\mu}\, H_{\kappa}\right.}
=\bigoplus_{i=0}^\infty\,\overline{b_   
{i}'}\ \overline{\mathbb{Z}}\]
  to a direct sum decomposition  
\[ 
K_\mu \left/ K_\mu \cap \bigcup_{\kappa<\mu}\,H_\kappa\right.=
\bigoplus_{{b_{\alpha_i}}\in
K_\mu\setminus\bigcup_{\kappa<\mu}\,H_\kappa}\, y_{i}'\,\mathbb{Z} 
\]
 with units $\left\{u_i\right\}$ such that $y_i'-u_i\,b_{i}'\in p^\nu
G$. Now set $y_{\alpha_i} = y_i' + u_i\, c_i$ so that $y_{\alpha_i}$
lifts $b_{\alpha_i}$.   

Assume for all $\mu<\lambda$, $H_\mu=\bigoplus_{\kappa\le\mu} L_\kappa$, where $L_\kappa$
is the free group generated by a lifting of the decomposition of $H_\kappa\left/\bigcup_{\kappa'<
\kappa}\, H_{\kappa'}\right.$ generated by the appropriate subset of $\mathfrak{B}$.  Then we  have
$\bigcup_{\mu<\lambda}\,H_\mu=\bigoplus_{\mu<\lambda}\,L_\mu$ and by the above, 
$H_\lambda=L_\lambda \oplus \bigoplus_{\mu<\lambda}\,L_\mu$.  By transfinite induction
we get $G=\bigoplus_{\lambda<\Omega}\,L_\lambda$.
\end{proof}

\subsection{Infinite Gaussian elimination modulo $p^\nu$.}

The reader is assumed thoroughly familiar with the details of Gaussian elimination as
developed in an introductory linear algebra course. Infinite Gaussian
elimination on a row finite $\omega \times \omega $ matrix can proceed very
much like the algorithm on a finite matrix. As in \cite{HoffmanKunze}, one
looks for a pivot in a row rather than a column as in many texts and standard implementations of finite Gaussian elimination. That insures that
only a finite number of entries need to be examined to either obtain a unit
pivot or to know that no such pivot exists. Subtracting multiples of a pivot
row from all other rows to make entries in the pivot column equal to 0 will, in general,
involve an infinite number of operations before the algorithm is complete.  To
avoid this, in the infinite case, 
rows are included with previously obtained pivot rows one at a time, 
and one clears the previously
obtained pivot columns in a row at the time that the row is included, and then 
finds a pivot if possible and clears above the pivot in the
new pivot column. In the infinite case there is no LU decomposition or forward pass and back
substitution because these might lead to rows changing infinitely often, and
there are no row permutations because some row might conceivably be permuted to a higher
numbered position an infinite number of times and thus never examined for a pivot. 
However, it is still the case
that a row finite $\omega \times \omega $ matrix is invertible if and only
if with these modifications of standard Gaussian elimination, infinite Gaussian
elimination will row reduce the matrix to a matrix whose columns are a permutation
of the columns of the identity matrix. 

We now modify infinite Gaussian elimination to produce
 an algorithm which we call infinite Gaussian elimination modulo $p^\nu$.\footnote{The
author has a working Maple V implementation of this algorithm. See the appendix in the
copy of this paper archived on http://arXiv.org or URL %\hfill\break
http://www.math.rutgers.edu/pub/osofsky/getbasis.html} This algorithm 
clearly also works if we have a finite
matrix $A$. We indicate the variables needed in the algorithm with a little information about them, 
 then give the steps of the algorithm, and then add a step by step 
explanation of what unusual steps do.   We start with a row finite  
$\omega\times\omega$ matrix  $\mathbf{A}$ with entries
in $\mathbb{Z}$.  In our proof of Theorem~\ref{mainthm}, the rows of
$\mathbf{A}$ will be some lifting of a given basis for 
$\overline{\mathbb{Z}}^{(\omega)}$ to elements of $\mathbb{Z}^{(\omega)}$.

By the expression
`principal submatrix' of an infinite matrix, we will
mean the submatrix obtained by taking the first $n$ rows and first $k$ columns
of the matrix, where $n$ and $k$ are both finite.  A `principal minor' will be 
the determinant of a square principal submatrix. 

Additional variables are needed to perform the algorithm. 
We use a diagonal matrix $\mathbf{U}$  
(or a countable row vector) to hold units modulo $p^\nu$.  Multiplying  row $i$
of $\mathbf{A}$ by an
appropriate unit $\mathbf{U}_{i,i}$ enables us to make a crucial determinant 1. 
The actual row reduction is done
in arbitrarily large but finite principal submatrices of an
$\omega\times \omega$ matrix $\mathbf{R}$.  Another 
$\omega\times \omega$ matrix $\mathbf{C}$ (for candidates) holds, in a finite
principal submatrix, the current  
candidates for lifting basis elements times units. These candidates change during
the elimination but each row only changes a finite number of times.  As the algorithm
progresses, we   
 multiply (an initial segment of) row $i$ of $\mathbf{A}$ by the appropriate unit $\mathbf{U}_{i,i}$ (integer relatively prime to
$p^\nu$) and then insert it into both $\mathbf{R}$ and $\mathbf{C}$.  All changes 
to $\mathbf{C}$ other than the concatenation of rows from $\mathbf{U}\mathbf{A}$ consist of adding multiples of $p^\nu$ to entries so nothing changes modulo $p^\nu$.   
In addition, we use   a finite
square matrix $\mathbf{M}$ which is generated from a submatrix of $\mathbf{C}$  
and has determinant 1. 

At the end of each loop of this algorithm, the matrix $\mathbf{R}$ will be a row reduction of $\mathbf{C}$
with row operations captured by $\mathbf{M}$.  Also, any entry of $\mathbf R$ which is a multiple
of $p^\nu$ is 0; it is set to 0 before any arithmetic is done using it.  
At any given stage of the algorithm
we work with  finite matrices large enough to hold all nonzero entries 
in a finite number of rows.  Moreover, the results of each loop of the algorithm applied to
$\overline{\mathbf{A}}$ are identical with the results of applying normal infinite Gaussian
elimination to $\mathbf{A}$.

\begin{algorithm}[Infinite Gaussian elimination modulo $p^\nu$] We start with an 
$\omega\times\omega$ integer valued row finite matrix $\mathbf{A}$. 

\begin{steplist}
\item Initialize.  Let your row index I be set to 0. Set up the matrix variables
$\mathbf{M}$, $\mathbf{C}$,  $\mathbf{R}$ and $\mathbf{U}$. Set up a row vector J to hold pivot columns. Read the $0^{th}$
row of $\mathbf{A}$ into $\mathbf{R}$, replacing any element divisible by $p^{\nu}$ with 0.

\item \label{new1} For K going from 0 to I -- 1, subtract $\mathbf{R}_{I,\,J(K)}$
times row K of $\mathbf{R} $ from row I of $\mathbf{R}$.

\item  \label{pass1} Search row I of $\mathbf{R}$ for the first entry which is relatively
prime to $p$.  If no such element is found
then STOP. The rows of $\mathbf{A}$ do not form a basis modulo $p^\nu$.
 Otherwise, let the first entry relatively prime to 
$p$ be in column J(I), and call column J(I)
the I$^{th}$ pivot column.

\item  Set $\mathbf{U}_{I,\,I}$ equal to an integer $u$ such that 
$u\mathbf{R}_{I,J(I)}\equiv1\mod p^\nu$.  Multiply row I of $\mathbf{A}$ by $u$.  If some entry in
the resulting row is a multiple of $p^\nu$, set that entry to 0. Insert the
result as row I in both $\mathbf{C}$ and 
$\mathbf{R}$.

\item For K going from 0 to I -- 1, subtract $\mathbf{R}_{I,\,J(K)}$
times row K of $\mathbf{R} $ from row I of $\mathbf{R}$.

\item \label{pivot1} The pivot in row I of $\mathbf{R}$ is now congruent to 1 modulo $p^{\nu}$.  Subtract a multiple of
$p^\nu$ from it to make the pivot 1.  Subtract the same multiple of $p^\nu$ from the
$\left({I,\,J(I)}\right)$ entry of $\mathbf{C}$.

\item  If any entry in row I of $\mathbf{R}$ is a multiple of $p^\nu$, subtract 
that multiple of $p^\nu$ from the corresponding entry in $\mathbf{C}$ and  set the
entry in $\mathbf{R}$ equal 0.

\item For K going from 0 to I -- 1, subtract $\mathbf{R}_{K,\,J(I)}$
times row I of $\mathbf{R} $ from row K of $\mathbf{R}$ to clear every entry in column
J(I) above the I$^{th}$ row.

\item \label{set0} If any entry in $\mathbf{R}$ is a multiple of $p^\nu$, then set that entry equal to 0.

\item \label{TheKey} Set $\mathbf{M}$ equal to the matrix $\left[\mathbf{C}_{\mathrm{K,\,J(K)}}
\right]_{0\le\mathrm{K}\le\mathrm{I}}$.  Set $\mathbf{C}=\mathbf{MR}$.

%%%%% Begin added step%%%%%
\item \label{remzero}
For each nonpivot column $\ell$ of $\mathbf{C}$, check to see if the
first nonzero entry $\mathbf{C}_{k,\ell}$ is divisible by $p^\nu$. 
If so, form the set 
$\mathfrak{S}_{\ell}$ consisting of all $l_i$ such that 
column $l_i$ is a pivot column, $\mathbf{C}_{k,l_i}$ is the first nonzero entry
in column $l_i$, and $\mathbf{R}_{l_i,\ell}\ne0$. 
If $\mathfrak{S}_\ell\ne\emptyset$, check if $p^\nu$ times the gcd of
$\mathfrak{S}_\ell$ divides $\mathbf{C}_{k,\ell}$.  If so,
express this gcd as a sum
$\sum_{\mathfrak{S}_{\ell}}\,\mathbf{C}_{k,l_i} b_{l_i}$.  Form a
column vector with zeros everywhere except for $b_{l_i}\cdot\left.
\mathbf{C}_{k,\ell}\right/d$ in row $l_i$, and add this to column $\ell$ of
$\mathbf{R}$. Premultiply by $\mathbf{M}$, and use the result as the new
column $\ell$ of $\mathbf{C}$.  The new $\mathbf{C}_{k,\ell}$ will be
0.
%%%%%End added step%%%%%

\item Read row I + 1 of $\mathbf{A}$ into $\mathbf{R}$, replacing multiples of $p^\nu$ by 0.   

\item \label{inc1} Increment I by 1 and GOTO \ref{new1}.
\end{steplist}
END \end{algorithm} 
\nobreak \enlargethispage{\baselineskip}

 That is the end of the algorithm.  To get a picture of what is happening, at the end of 
the $(n-1)^{th}$ loop at {\ref{inc1}} the column permuted matrix  $\mathbf R$ 
(picturing $j(i)$ as
though it were $i$) looks like
\[ \begin{array}{c} $\huge R$ \end{array} = 
\left[  \begin{array}{*{20}c}
   \begin{array}{*{20}c}
   \ 1\  & \ 0\  & \cdots  & \ 0\  &\vline &  {r_{0,n} }\vphantom{\Bigg)}  \\
   \ 0\  & \ 1\  & \cdots  & \ 0\  &\vline &  {r_{1,n} }  \vphantom{\Bigg)} \\
    \ \vdots\   &  \ \vdots\   &  \ \ddots\   &  \ \vdots\   &\vline &   \vdots \vphantom{\Bigg)}  \\
    \ 0\  &  \ 0\ & \cdots  & \ 1\  &\vline &  {r_{n - 1,n} } \vphantom{\Bigg)} \\
\hline
   {r_{n,0}} & {r_{n,1}} &  \cdots  & {r_{n,n-1} } &\vline &  {r_{n,n} } \vphantom{\Bigg)} \\
\end{array} &\vline &  \begin{array}{*{20}c}
   {\vphantom{a_1}}  \vphantom{\Bigg)} \\
   \ \ $\huge{B}${\vphantom{m_1}} \qquad \cdots \vphantom{\Bigg)} \\
   {\vphantom{m_1}} \vphantom{\Bigg)} \\
%   {\vphantom{m}}    \\
\hline
   $\bf D${\vphantom{m_1}} \qquad \cdots \vphantom{\Bigg)} \\
\end{array}  \\
\end{array} \right]
\]
for an appropriate finite matrix $\mathbf{B}$ and finite row $\mathbf{D}$, and all
entries in $\mathbf{R}$ which are divisible by $p^{\nu}$ are 0. \medskip

Now for a more detailed explanation of how this algorithm works. In the permuted matrix used in the discussion, $j(i)$ will be treated as though it were $i$ to aid in visualization
of the progress of the algorithm. That is, we will pretend that we have permuted the
columns of the matrix. \medskip

\ref{new1} is the first pass at clearing already obtained pivot columns (which have
pivot 1) in row $i$.  It is used to get the unit mod $p^\nu$ we must multiply the $i^{th}$
row of $\mathbf{A}$ by to make sure that we can make the pivot in row $i$ equal to 1. It is not performed when $i=0$. 
\medskip

\ref{pivot1} relies on the claim that the pivot is congruent to 1 modulo $p^\nu$. Why is 
that claim true?  Adding one row of a matrix to another corresponds to premultiplication
by a matrix of determinant 1.  After \ref{pass1}, if we look at the principal minor
of the column permuted matrix $\mathbf{R}$, it has determinant the $\left(i,\,i
\right)$ entry of the permuted
$\mathbf{R}$ because it is upper triangular with all other diagonal entries 1.  
When we multiply what was the last row before \ref{pass1} by $u$, we make that determinant congruent to 1 modulo
$p^\nu$.  Now we redo the elementary row operations of determinant 1 to get an upper 
triangular  matrix with element in the $\left(i,\,j(i)\right)$ slot equal to the determinant.\medskip

In \ref{pivot1}, subtracting multiples of $p^\nu$ from the same entries in both $\mathbf{R}$
and $\mathbf{C}$ does not change  $\overline{\mathbf{C}}$ and does insure that the 
elementary row operations we have done so far will reduce the new $\mathbf{C}$ to the
new $\mathbf{R}$.\medskip

Since we want entries in $\mathbf{R}$ congruent to 0 mod $p^\nu$ to be 0, we set them to 0 in \ref{set0}.  This can only affect entries in nonpivot columns.  Now we must make sure that
our $\mathbf{C}$ row reduces to the new $\mathbf{R}$.  This is done in \ref{TheKey}.  At this stage, the appropriate principal submatrix of the column permuted matrix $\mathbf{R}$ 
is the identity matrix.  So the row operations we have done have reduced the 
corresponding principal submatrix of the column permuted matrix $\mathbf{C}$ to the identity.  
By standard linear algebra,  the matrix $\mathbf{M}$ is the inverse of the product of the
elementary matrices which produce this elimination by premultiplication.  Thus from
$\mathbf{R}=\mathbf{M}^{-1}\mathbf{MR}$ we see that setting $\mathbf{C}=\mathbf{MR}$ 
gives us a matrix which row reduces to the new $\mathbf{R}$, and since $\mathbf{R}$ did not
change modulo $p^\nu$, neither did $\mathbf{MR}$.

%%%%%Begin added explanation%%%%%
In \ref{remzero}, the algorithm bounds the power of $p^\nu$ that can divide entries of 
$\mathbf{C}$ after the corresponding row of $\mathbf{A}$ becomes all zeros. This step may
change $\mathbf{C}$ and nonzero entries in $\mathbf{R}$ modulo $p^\nu$.  
If the first entry $\mathbf{C}_{k,l}$ in a nonpivot 
column of $\mathbf{C}$. The several imposed conditions on $\mathfrak{S}_\ell$ 
insure that no zero entry of $\mathbf{R}$ becomes nonzero, and the divisibility 
property makes the added vector a multiple of $p^\nu$.  If a nonzero entry appears
in $\mathbf{C}$ after all the nonzero mod $p^\nu$ entries in its row occur in pivot columns,
it may propagate, but that leads to entries in the row divisible by higher powers of
$p^\nu$, and eventually \ref{remzero} will make all of these entries zero. Thus \ref{remzero} makes sure that no row has an infinite number of entries congruent to 0 modulo $p^\nu$.

New row operations are only done to the rows above the pivot row when their entries in the current pivot column is nonzero.   Hence
once the finite set of rows of $\mathbf{R}$ from $0$ to $i$ have zero
entries except for a pivot of 1, 
%%%%%additional part of the explanation%%%%%
and there are no more nonzero multiples of $p^\nu$ in these rows of $\mathbf{C}$,
%%%%%end additional part of the explanation%%%%%
those rows will no longer be affected by the elimination process.\medskip

The last steps of the algorithm just set up for the next loop.  
\bigskip

\subsection{ The proof of Theorem~\ref{mainthm}.}

\begin{proof} By Corollary~\ref{cgenuf}, it is enough to show that, for a countably generated
free abelian group $G$ with $\mathfrak
{B}$ a basis for $\overline{G}$,  there is a direct decomposition lifting of
\[ \overline{G}=\bigoplus_{\alpha\in\mathfrak{B}}\,\overline{b_\alpha}\,\overline{\mathbb{Z}}
\]
to the direct decomposition
\[ G= \bigoplus_{\alpha\in\mathfrak{B}}\,y_\alpha\mathbb{Z}
\]

Form a matrix $\mathbf{A}$ whose rows are some lifting of
$\mathfrak{B}$.
   Do infinite Gaussian elimination modulo $p^\nu$ on 
$\mathbf{A}$.  Since the rows of $\overline{\mathbf{A}}$
form a basis for $\overline{\mathbb{Z}}^{(\omega)}$ and modulo $p^\nu$ this algorithm agrees
with infinite Gaussian elimination, after a finite number of steps, the top $i+1$ rows
of $\mathbf{R}$ will be rows of the identity and all rows of the identity will eventually
arise
as rows of $\mathbf{R}$. Since every entry of $R$ which is zero modulo $p^\nu$ is actually $0$, $\mathbf{C}$ is row reduced to the identity provided every row
at some point stops changing in taking the product $\mathbf{MR}$. 
Since all of the entries of row $n$ of $\mathbf{C}$ which are not congruent to $0$ mod $
p^\nu$ are contained in a finite number of columns, 
 any row of $\mathbf{C}$ ceases to change when all the rows of the identity with $1$ in 
those columns have been obtained in the matrix $\mathbf{R}$.
Hence after an infinite number of steps each row of $\mathbf{C}$ will have stabilized and
the stabilized rows of $\mathbf{C}$ will form a basis for $\mathbf{Z}^{(\omega)}$ which lifts the direct sum decomposition.
\end{proof}

\section{Lattices of commuting idempotents}\label{Part3}

\subsection{Definitions and notation.}
The following notation will be used, usually without comment, in the rest of this paper.

Let $\mathcal{E}$ be a lattice of commuting idempotents in a ring $R$ with $1$, that is, $\mathcal{E}$ is closed under multiplication and
addition of orthogonal idempotents. The idempotents in $\mathcal{E}$
together with the identity generate a boolean algebra $\mathcal{B}$ under
multiplication as in $R$ but addition the symmetric difference $e +_{
\mathcal{B}} f= e\left(1-f\right) + f\left(1 - e\right)$. Let $\mathbb{Z}
\left[ \mathcal{B}\right] $ be the semigroup algebra of $\left\langle 
\mathcal{B},\,\cdot \right\rangle $, that is, the free abelian group with
basis the elements of $\mathcal{B}$ and multiplication the multiplication as
in $\mathcal{B}$. Let 
\[
\mathcal{S}=\mathbb{Z}\left[ \mathcal{B}\right] \left/ \left\langle \left(
e+f\right) -e-f:ef=0\right\rangle \right.. 
\]
$\mathcal{S}$ is a free lattice ring in the sense that it can be formed for
any modular, complemented lattice and has appropriate universal properties
with respect to embedding such lattices in rings.

For convenience, we will assume that $\mathcal{E}$ is a Boolean ideal, that is, if $f=f^2\in\mathcal{E}
R$, then $f\in\mathcal{E}$.  This does not change $\mathcal{E}
R$. \medskip

\subsection{Elementary properties of $\mathcal{S}$.}

Much of the known material assumed
in this subsection can be found in graduate level text books such as \cite{lam}.

The next proposition is essentially a sequence of remarks, included with
short proofs.

\begin{proposition} \label{mainprop}
The following hold for the free lattice ring $\mathcal{S}$.

\begin{enumerate} \renewcommand{\theenumi}{(\alph{enumi})} \renewcommand{\labelenumi}{\theenumi}

\item The additive group of $\mathcal{S}$ is torsionfree.
 
\item  The lattice of idempotent generated ideals of $\mathcal{S}$ is
isomorphic to $\mathcal{B}$.

\item\label{projideal}  Any finitely generated ideal of $\mathcal{S}$ is 
cyclic and isomorphic to a sum $\sum_{i=1}^n \,f_i\mathcal{S}$ for some set
of orthogonal idempotents $\left\{f_i\right\}\subseteq\mathcal{B}$. 

\item  $R$ is an $\mathcal{S}$-module under the map induced by the inclusion
of $\mathcal{B}$ in $R$.

\item  The projective dimension of an idempotent generated ideal $I$ of $
\mathcal{S}$ is greater than or equal to the projective dimension over $R$
of the module $I\otimes _{\mathcal{S}}R$.
\end{enumerate}
\end{proposition}

\proof
\begin{enumerate} \renewcommand{\theenumi}{(\alph{enumi})} \renewcommand{\labelenumi}{\theenumi}
\item The kernel of the ring map from $\mathbb{Z}$ to $\mathcal{S}$ is
generated by idempotents and so pure.

\item  Any element of $\mathcal{S}$ is of the form $\sum_{i=1}^ne_in_i$
where $\left\{ e_i\right\} \subseteq \mathcal{B}$ are pairwise orthogonal
and $n_i\in \mathbb{Z}$. Assume such an element is idempotent. By the
torsionfree property of $\left\langle \mathcal{S},\,+\right\rangle $, the $
n_i$ must be all $1$, and $e=\sum_{i=1}^ne_i\in \mathcal{E}$. But then the
symmetric difference of $e$ and $f$ is the same as in $\mathcal{B}$.

\item  Given a finite set of idempotents $\left\{ e_i:1\le i\le n\right\}
\subseteq {\mathcal{B}}$, the minimal nonzero idempotents in the lattice
they generate will be pairwise orthogonal and generate the same lattice.
Since $\mathcal{S}$ is a quotient of the ring $\mathbb{Z}\left[\mathcal{B}\right]$,
any element of $\mathcal{S}$ is of the form $\sum_{j=1}^m e_jn_j$.  Moreover, if the
$\left\{e_j\right\}$ happen to be orthogonal,  $\left(\sum_{j=1}^m e_jn_j\right)
\mathcal{S}=\sum_{j=1}^m\,\left(e_jn_j\mathcal{S}\right)$. 

Now let $I$ be the finitely generated ideal 
\[I=\sum_{i=1}^k\left(\sum_{j=1}^{l_j}
\,e_{i,j}n_{i,j}\mathcal{S}\right)\subseteq \mathcal{S}.\]
Split each $e_{i,j}$ into an orthogonal sum of the nonzero
minimal elements in the lattice generated by $\left\{e_{i,j}: 1\le j\le l_j,\, 
1\le i \le k\right\}$.  Collecting multiples of each of these minimal elements,
we get a generator for $I$ of the form $\sum_{i=1}^{k^{\prime
}}f_im_i$ where the $\left\{ f_i\right\} $ are pairwise orthogonal
idempotents in $\mathcal{E}$. But $\left( \sum_{i=1}^{k''}f_im_i\right) 
\mathcal{S}\approx \bigoplus_{i=1}^{k^{\prime \prime }}f_i\mathcal{S}$ if we ignore
terms with $m_i=0$.

\item  The obvious map $\mathbb{Z}\left[ \mathcal{B}\right] \longrightarrow R
$ is a ring homomorphism whose kernel contains 
\[
\left\langle \left( e+f\right) -e-f:ef=0\right\rangle .
\]

\item  $I$ is a direct limit of idempotent generated cyclics and so flat. A
projective resolution 
\[
\cdots \rightarrow P_i\rightarrow P_{i-1}\rightarrow \cdots \rightarrow
P_0\rightarrow I\rightarrow 0
\]
is therefore pure exact. Moreover, since $P_i$ is a projective $\mathcal{S}$-module,
$P_i\otimes _{\mathcal{S}}R$ is a projective $R$-module. Thus 
\[
\cdots \rightarrow P_i\otimes _{\mathcal{S}}R\rightarrow P_{i-1}\otimes _{
\mathcal{S}}R\rightarrow \cdots \rightarrow P_0\otimes _{\mathcal{S}
}R\rightarrow I\otimes _{\mathcal{S}}R\rightarrow 0
\]
is a projective resolution of $I\otimes _{\mathcal{S}}R$.  If the kernel of 
a map $P_i\rightarrow P_{i-1}$ is $\mathcal{S}$-projective, by pure 
exactness and the fact that tensoring preserves projectivity we see that
the kernel of  
$P_i \otimes_{\mathcal{S}}R\rightarrow P_{i-1} \otimes_{\mathcal{S}}R$ is $R$-projective.  Thus the $\mathcal{S}$-projective dimension
of $I$ is at most $i$ implies that the $R$-projective dimension of $I\otimes_{\mathcal{S}}R$
is also at most $i$.  \endproof

\end{enumerate}
\medskip

\begin{proposition} \label{freegp} The
additive group of $\mathcal{S}$ is a free abelian group.
\end{proposition}

\proof   
Let $\mathfrak{X}$ be the family of all subsets $X$ of $\mathcal{B}\setminus
\left\{ 0\right\} $ such that whenever $\left\{ e_{i}\right\} $ is a set of
orthogonal idempotents in $X$, if $\left\{ f_{j}\right\} $ is any set of
orthogonal idempotents such that $\left\{ e_{i}\right\} \ne \left\{
f_{j}\right\} $ and $\sum_{i}e_{i}=\sum_{j}f_{j}$, then at least one $%
f_{j}\notin X$. $\mathfrak{X}$ is an inductive poset under $\subseteq $, so
by Zorn's lemma there is a maximal element $B$ in $\mathfrak{X}$.  $B$ is $%
\mathbb{Z}$-linearly independent in $\mathcal{S}$ because the only relations on
the $\mathbb{Z}$-linearly independent 
idempotents in $\mathbb{Z}\left[ \mathcal{B}\right] $ set an idempotent equal to
an orthogonal sum of other idempotents. $B$ will be a vector space basis for $%
\mathcal{B}$ over the field of 2 elements. Let $f\in \mathcal{B}\setminus
\left\{ 0\right\} $. If $f\notin B$, then $B\cup \left\{ f\right\} \notin %
\mathfrak{X}$. Hence there must be a set $\left\{ e_{i}\right\} $ of
orthogonal idempotents in $B\cup \left\{ f\right\} $  and a different
 set $\left\{ f_{i}\right\} \subseteq B\cup \left\{ f\right\} $ of
orthogonal idempotents with $\sum_{i=1}^{n}e_{i}=\sum_{j=1}^{m}f_{j}$. If $%
f\in \left\{ e_{i}\right\} \cap \left\{ f_{j}\right\} $ then we get $%
\sum_{e_{i}\ne f}e_{i}=\sum_{f_{j}\ne f}f_{j}$ with all summands in $B$, a
contradiction. Similarly, if $f\notin \left\{ e_{i}\right\} \cup \left\{
f_{j}\right\} $ we get a contradiction. Hence $f$ is in precisely one of the
two sets, say $f=e_{1}$. Then $f=\sum_{j}f_{j}-\sum_{i=2}^{n}e_{i}$ is in
the span of $B$. 
\endproof
\medskip

Proposition~\ref{freegp} strongly reinforces the observation that $\mathcal{S}$ is a
free object.  The basis found for its additive group will be a basis for $\mathcal{S}\otimes
_\mathcal{S}F$ over $F$ for any field $F$.

In his proof in  of the affirmative answer to the Wiegand question in the
case $n=1$, R. S. Pierce proved the next lemma with completely different
terminology. See \cite[Lemma 2.7]{pierce:76}.

\begin{proposition} \label{indep}
Let $\left\{ \kappa _{\alpha }\right\} $ be a set of elements in a submodule
of a free $\mathcal{S}$-module $K$, where the $\left\{ \kappa _{\alpha }\otimes 1\right\} $ 
are all nonzero. Then if $\left\{ \kappa _{\alpha
}\otimes _{\mathcal{S}}R\right\} $ is $R$-independent in $K\otimes _{
\mathcal{S}}R$, then $\left\{ \kappa _{\alpha }\right\} $ is 
$\mathcal{S}$-independent in $K$.
\end{proposition}

\proof
Assume not. Then there is a shortest sum $\sum_{i=1}^{n}\,\kappa _{\alpha
_{i}}s_{i}=0$ where the summands are all nonzero in $\mathcal(S)$. Considering elements of the free module $K$ as consisting of
sums of idempotents times basis elements, we see that the annihilator of
each $\kappa _{\alpha _{i}}s_{i}$ is generated by an idempotent $\left(
1-\varepsilon _{i}\right) $. Since $n$ is the smallest number of summands
that can give you a zero and $\sum_{i=1}^{n}\,\kappa _{\alpha
_{i}}s_{i}\varepsilon _{1}=0$, we have $\,\kappa _{\alpha
_{i}}s_{i}\varepsilon _{1}\ne 0$ for all $i$. Similarly $\kappa _{\alpha
_{i}}s_{i}\varepsilon _{1}\varepsilon _{2}\ne 0$ for all $i$. Continuing in
this manner we get $\kappa _{\alpha _{i}}s_{i}\prod_{j=1}^{n}\varepsilon
_{j}\ne 0$ for all  $i$. Then $\sum_i\,\kappa_i s_i \prod_{j=1}^{n}\varepsilon _{j}$ has all
summand nonzero and there is an integer $m$ such that $
\sum_{i=1}^{n}\,\kappa _{\alpha _{i}}s_{i}m^{-1}\prod_{i=1}^{n}\varepsilon
_{i}$ is an element not divisible by any integers other than $\pm 1$ in the free abelian
additive group of $K$. But then $\sum_{i=1}^{n}\,\kappa _{\alpha
_{i}}s_{i}m^{-1}\prod_{i=1}^{n}\varepsilon _{i}\otimes 1$ is nonzero in $%
K\otimes _{\mathcal{S}}R$ and each of the summands is nonzero. \endproof

We quote a Proposition due to Kaplansky that is basic to almost all
studies of infinitely generated projective modules, with two consequences
giving rise to the same result for von Neumann regular rings. 

\begin{proposition}[Kaplansky]
\label{quote}A projective module over any ring is a direct sum of countably
generated submodules. From this we obtain:

\begin{enumerate}
\renewcommand{\theenumi}{(\alph{enumi})} \renewcommand{\labelenumi}{%
\theenumi}

\item  Any projective right module over a von Neumann regular ring is
isomorphic to a direct sum of cyclic (idempotent generated) right ideals.

\item  Any projective module over a commutative semihereditary ring is
isomorphic to a direct sum of finitely generated right ideals.
\end{enumerate}
\end{proposition}

See \cite{Kaplansky} for a proof. The proof of this theorem is the template
on which the preliminary proofs in Section~\ref{kap} are based.

\subsection{The proof of an affirmative answer to the Wiegand question.}

We now complete our work on the Wiegand question.  

\begin{proposition}
Let $R$ be a commutative von Neumann regular ring. Let $F$ be a projective $%
\mathcal{S}$-module and let $K$ be any pure submodule of $F$. Then if $K\otimes _{%
\mathcal{S}}R$ is projective as an $R$-module, then $K$ is projective as an $%
\mathcal{S}$-module.
\end{proposition}

\proof  Since $K\otimes _{\mathcal{S}}R$ is a projective $R$-module, it is a direct
sum of the form $K\otimes _{\mathcal{S}}R=\bigoplus_{\alpha }\thinspace x
_{\alpha }R$ where for each $\alpha $ there is an $e_{\alpha }$ such that $%
x_{\alpha }R\approx e_{\alpha }R$. If any $e_{\alpha }$ is of finite but
composite order, express it as an orthogonal sum of idempotents of prime
power order by the Chinese Remainder Theorem. In the von Neumann regular
case where there are no nilpotent elements, the prime power must be the
prime itself. We can then divide the indexing set into a family of subsets

\[
\mathfrak{F}_{p}=\left\{ \alpha :\mathop{char}\left( e_{\alpha }\otimes _{%
\mathcal{S}}R\right) =p\right\} 
\]
for $p$ a prime or $0$.

Consider the map $K\stackrel{I_{k}\otimes 1}{\longrightarrow }K\otimes _{%
\mathcal{S}}R\longrightarrow \bigoplus_{\alpha \in \mathfrak{F}%
_{0}}\thinspace x_{\alpha }R$. Its image is a projective $\mathcal{S}$%
-module, so it splits. Hence without loss of generality we can work with the
kernel of this map in place of $K$ and assume that $K\otimes _{\mathcal{S}}R$
is torsion. But then it is the orthogonal sum of its $p$-primary components
so we need only look at sums of the form $\bigoplus_{\alpha \in \mathfrak{F}%
_{p}}\thinspace x_{\alpha }R$ for a fixed prime $p$. That is, without
loss of generality, $K\otimes _{\mathcal{S}}R$ is $p$-primary. Since the
additive group of $\mathcal{S}$ is free, the additive group of $F$ is free
and hence $K$ is a subgroup of a free abelian group and so free.  By Theorem~%
\ref{mainthm}, there is a basis $\left\{ b_{\lambda }\right\} $ of $K$ which
lifts the direct sum decomposition $G_{p}\left/ pG_{p}\right. =
\bigoplus_{\alpha \in \mathfrak{F}_{p}}\thinspace x_{\alpha }R$ to a
direct sum decomposition of $K$. 

For every $\alpha $, let  $\mathfrak{B}_{\alpha }=\left\{ b_{\lambda
}:b_{\lambda }\otimes 1\in x_{\alpha }R\right\} $. 
Let $H_{\alpha }$ be the $\mathcal{S}$-submodule of $K$ generated by $\mathfrak{B}
_{\alpha }$.  Since the generators of $H_{\alpha }$ all map to $x_{\alpha }R$
under $Id_{K}\otimes 1_{R}$, so must $H_{\alpha }$. 
Since $H_\alpha $ contains $\mathfrak{B}_\alpha$ and
$\bigcup_\alpha\mathfrak{B}_\alpha$ is a basis for
$K$, $K=\sum_\alpha\, H_\alpha$. By Proposition~\ref{indep}, that sum
is direct.

Select any element $y$
in $H_{\alpha }$ which maps to $x_{\alpha }$. This $y$ is an element lying
in a finitely generated free submodule of $F$.  Hence it is of the form $%
y=\sum_{i=1}^{m}\sum_{j=1}^{k_{i}}c_{i,j}e_{i,j}n_{i,j}$ where the $c_{i,j}$
are basis elements of $F$, and we can use our little trick of decomposing
into the minimal idempotents in a finite lattice to get that $e_{i,j}$ and $%
e_{k,l}$ are either the same idempotent or orthogonal. Because of the $\mathbb{Z%
}$-purity of $K$, we may find a $y_{\alpha }\in H_{\alpha }$ such that each
sum of the form $\sum_{e_{i,j}=e_{k,l}}c_{i,j}e_{i,j}n_{i,j}$ is of content
1 and hence this  $y_{\alpha }$ generates a direct summand of $F$. But
then $y_\alpha\,\mathcal{S}$ is a direct summand of $H_\alpha$ which
maps to the same submodule of $K\otimes_\mathcal{S}R$.  We conclude
that $H_\alpha=y_\alpha\,S$ for all $\alpha$. 
Thus $K=\bigoplus _{\alpha }y_{\alpha }\mathcal{S}$ so $K$ is projective.
\endproof

\begin{corollary}\label{projres}
Let $F$ be a projective $\mathcal{S}$-module of the form 
\[ F = \bigoplus_{\alpha\in\mathfrak{A}}\, e_\alpha\mathcal{S} \]
where each $e_\alpha\mathcal{S}$ is isomorphic to an ideal of $\mathcal{S}$ contained
in $\mathcal{ES}$. Then for any pure submodule $K$ of $F$, 
$\mathop{\rm pd}\nolimits_R\left( K\otimes _{\mathcal{S}}R\right) =
\mathop{\mathrm{pd}}_{\mathcal{S}}\left( K\right) $.
\end{corollary}

\begin{proof}  We can take a short projective resolution of $K$ over $\mathcal{S}$, say 
\[
0\longrightarrow L\longrightarrow P\longrightarrow K\longrightarrow 0 
\]
is exact with $P$ projective and, like $F$, a direct sum of cyclic projectives
of the form $e\mathcal{S}$ for some $e\in\mathcal{E}$. Then if we let $\infty -1=\infty $, $
\mathop{\rm pd}\nolimits_{\mathcal{S}}\left( L\right) =\mathop{\rm pd}
\nolimits_{\mathcal{S}}\left( K\right) -1$. This short exact sequence is
pure, so tensoring with $R$ over $\mathcal{S}$ gives a short projective
resolution of $K\otimes _{\mathcal{S}}R$
\[
0\longrightarrow L\otimes _{\mathcal{S}}R\longrightarrow P\otimes _{\mathcal{
\ S}}R\longrightarrow K\otimes _{\mathcal{S}}R\longrightarrow 0 
\]
with $\mathop{\rm pd}\nolimits_{R}\left( L\otimes _{\mathcal{S}}R\right) =
\mathop{\rm pd}\nolimits_{R}\left( K\otimes _{\mathcal{S}}R\right) -1$.
Induction on $\mathop{\rm pd}\nolimits_{\mathcal{S}}\left( K\right) $
completes the proof.
\end{proof}

\begin{theorem}[The answer to the Wiegand question]
\label{last}For any commutative von Neumann regular ring $R$ with a commuting set of idempotents $\mathcal{E}
$, $\mathop{\rm pd}\nolimits_R\left( \mathcal{E}R\right) =\mathop{\rm pd}
\nolimits_{\mathcal{S}}\left( \mathcal{E}\mathcal{S}\right) =\mathop{\rm pd}
\nolimits_{\mathcal{B}}\left( \mathcal{EB}\right) $.
\end{theorem}

\begin{proof} $\mathcal{ES} $ has a projective resolution of the form required in
Corollary~\ref{projres}.  Then Corollary~\ref{projres} gives the desired conclusion.
\end{proof}

One way to summarize this answer to the Wiegand question is to say that, when
working in a submodule of a free module over a commutative regular ring, the 
lattice of direct summands carries all of the information about the module, and
the coefficients essentially none. For example, 
note that in Theorem~\ref{last}, the lattices of direct summands in the
three ideals $\mathcal{E}R$, $\mathcal{E}\mathcal{S}$, and $\mathcal{EB}$
are isomorphic, as they correspond to the idempotents themselves. However,
as soon as one gets to free modules on more than one generator, that
property fails. Since the number of one dimensional subspaces of a
2-dimensional vector space depends on the
cardinality of the field, if $R=\mathcal{S}\left/3\mathcal{S}\right.$ then the number of direct summands of $eR\oplus eR$ isomorphic to $eR$
and the number of direct summands of $e\mathcal{B}\oplus e\mathcal{B}$ isomorphic to $e
\mathcal{B}$ will always be different for any idempotent $e$.

\medskip
\section{Appendix}
Here is a Maple program which implements an algorithm similar to but not identical
with the infinite gaussian elimination modulo $p^\nu$ of this paper.  
The $\mathbf{C}$ of this
algorithm is the analogue of the $M$ in the algorithm here. Except for pivot columns,
$\mathbf{R}$ is only determined modulo $p^\nu$. The \# indicates a comment in the program.
The program for the algorithm in this paper, as well as in this appendix, 
can be found via URL {\tt http://math.rutgers.edu/$\sim$osofsky} in both .mws and .html
formats.
\smallskip

\parindent=0pt
\hrulefill 

\smallskip
\#  mgcdex IS A PROGRAM TO COMPUTE THE GREATEST COMMON DIVISOR OF A

\#  VECTOR OF INTEGERS, AND A LINEAR COMBINATION OF ENTRIES OF THE 

\#  VECTOR WHICH GIVES THAT GCD.
\smallskip

mgcdex:=proc(A,B)
local i,j,a,b;
with(linalg):
\smallskip

B:=array(1..vectdim(A)+1);
for i from 1 to vectdim(A) do
   B[i]:=1; B[vectdim(A)+1]:=0;
od; \smallskip

for i from 1 to vectdim(A) do

\hspace{.25truein} B[vectdim(A)+1]:=igcdex(B[vectdim(A)+1],A[i],'a','b');
   B[i]:=b;  
 
\hspace{.25truein}   for j from 1 to i-1 do
      B[j]:=a*B[j];
   od; 
             
od;
       
end: \bigskip

\# {\bf THE PROGRAM  GetBasis IMPLEMENTING A VARIANT OF 

{\rm\#} GAUSSIAN ELIMINATION MOD $p^\nu$.}\bigskip

\# The input consists of a finite matrix A and a prime power $p$.\smallskip

GetBasis:=proc(A,p) 
  local
ind, i, j, checkdet, k, temp, n, m, h, l, mat, Adj, u,check, mat1, getgcd, ell, hold, 
ii, V, B, fl0, sum: 
  global R, C, U, mgcdex: with(linalg):\medskip

\# {\bf INITIALIZE}\medskip

\# For checking purposes we will also hold the inverse of A in R.

\# C is the matrix whose rows are the required basis. 

\# The second part of the augmented matrix B will be the 

\# inverse of C.  It is not necessary to do this but it may help. \smallskip
 
R:=array(1..rowdim(A), 1..coldim(A)):
copyinto(A,R,1,1):   C:=array(1..rowdim(A),1..rowdim(A)); 
 for i from 1 to rowdim(A) do for j from 1 to coldim(A) do
     R[i,j]:=mods(R[i,j],p):
 od:od: \smallskip

\# ind(ex) holds our column permutation. \smallskip

 ind:=array(1..coldim(A)): for i from 1 to coldim(A) do ind[i]:=i; od; \smallskip

\# U is a diagonal matrix of units modulo p used to multiply 

\# rows and make pivots 1. \smallskip

U:=array(1..rowdim(A),1..rowdim(A));     
 for i from 1 to rowdim(A) do for j from 1 to rowdim(A) do
         if (i<>j) then  C[i,j]:=0; U[i,j]:=0; else           
         C[i,j]:=1; U[i,j]:=1;   fi:       
  od; od; 
  R:=concat(R,C);  \smallskip

\# checkdet holds candidates for the next pivot. \smallskip

checkdet:=array(1..coldim(A)):        
  if (rowdim(A)>coldim(A)) then RETURN(`More rows than columns cannot
form a basis.`):fi: \smallskip

\# We need a temporary location to compute changes in C to avoid

\# nonzero entries divisible by p. \smallskip

  getgcd:=array(1..rowdim(A)); \smallskip

\# The index 'i' will stand for the row currently being worked on. \smallskip

\# {\bf END INITIALIZE}
\medskip

\# {\bf THE ACTUAL COMPUTATION} \smallskip

\# The variable $i$ will denote the working row.

\# Compute the determinant of the block to be used and the pivot

\# column by looking for a unit mod $p$ to be the next pivot. \smallskip

for i from 1 to rowdim(A) do 
   
      for m from 1 to coldim(A) do checkdet[m]:=R[i,m]:       
         for k from 1 to i-1 do 
             checkdet[m]:=checkdet[m]- R[i,ind[k]]*R[k,m]:       
         od:
      od:      
      for k from 1 to coldim(A) while (igcd(checkdet[k],p) $<>$ 1) do :
od:  

if (k $>$ coldim(A)) then print(A,R,C):RETURN(`No pivot. Not a basis
mod p.`): fi: \smallskip

\# If necessary, permute columns by permuting entries of ind. \smallskip

 if (ind[i]<>k) then for n from 1 to coldim(A) while (ind[n] $<>$ k) do : 
od:
 temp:=ind[i]:ind[i]:=k:ind[n]:=temp:  fi: \smallskip

\# Multiply the working row by the inverse of the pivot to 

\# make the pivot 1 mod p and clear below the diagonal. \smallskip

 u:= (checkdet[k]\^(-1) mod p) : U[i,i]:=
mods(u,p):   
      for n from 1 to coldim(A) do R[i,n]:= mods(u*R[i,n],p): od: \smallskip

\# Clear below the permuted diagonal. \smallskip

if (i<>1) then 
   for h from 1 to i-1 do         
        temp:=mods(R[i,ind[h]], p);
        R:=addrow(R,h, i, -temp): 
        C:=addcol(C, i, h, temp)
    od:
  fi:  
\smallskip

\# Clear above the diagonal. \smallskip

if (i $<>$ 1) then   
    for h from i-1 to 1 by -1 do 
       temp:=mods(R[h,ind[i]],p);
       R:=addrow(R,i,h,-temp): 
       C:=addcol(C,h, i, temp);
    od:  
  fi:  \smallskip

\# We now correct for some (enough) nonzero multiples of p which may

\# occur in our candidate C for a lifting. \smallskip

m:=0; ell:=0; fl0:=0;      

if (mods(C[k,ind[i]],p)$<>$ 0) then fl0:=1; fi; 
 
if ((C[k,ind[i]] $<>$ 0)and(igcd(C[k,ind[i]],p)=p)) then

\hspace{.25in} for j from 1 to i-1 do sum:=0;

\hspace{.25in} if (i $<>$ 1) then 
                      for ii from 1 to k-1 do
                         sum:=sum+C[ii,j]\^{}2;   
                      od;
                   fi;  
                   
\hspace{.25in} if (sum=0) then ell:=ell+1; getgcd[ell]:=ind[j];fi; od;

\hspace{.25in} if (ell $<>$ 0) then 

\hspace{.25in}\hspace{.25in}
            V:=array(1..ell); for j from 1 to ell do
V[j]:=C[k,getgcd[j]]; od;
            
\hspace{.25in}\hspace{.25in} mgcdex(V,B); 

\hspace{.25in}\hspace{.25in}
            if (mods(C[k,ind[i]],(B[ell+1]*p))=0) then m:=1; else
fl0:=1; fi;
         
\hspace{.25in}
fi;

\hspace{.25in} if (m=1) then  k:=k-1; temp:=C[k,ind[i]]/B[ell+1];  
    for j from 1 to ell do

\hspace{.25in}\hspace{.25in}        C:=addcol(C,getgcd[j],i,(-B[j]*temp));
        R:=addrow(R,i,getgcd[j],B[j]*temp);

\hspace{.25in}    od; 

fi; \smallskip

\# Printouts added to observe progress. \smallskip

mat:=submatrix(C,1..i,1..i);
mat1:=submatrix(R,1..i,1..coldim(A)+i); 
 print(`Row `,i,` C = `,mat,` Rowreduction = `,mat1):  \smallskip 
 
od:  \# This is the end of the working program. \smallskip

  print(`Orig A = `,A,` C = `,C,` U = `,U,`R = `,R);end: 

\hrulefill

\end{document}